\input ./style/arxiv-general.cfg
\documentclass[aap,MSNbibl,dvips]{arximspdf}
\makeatletter
   \@ifpackageloaded{graphicx}{}{\usepackage{graphicx}}
\makeatother
\usepackage{mathbh}

\doi{10.1214/14-AAP1095}
\volume{26}
\issue{1}
\pubyear{2016}
\firstpage{402}
\lastpage{424}
\docsubty{FLA}

\makeatletter
\renewcommand{\P}{\mathbb{P}}
\newcommand{\lleft}{\left}
\newcommand{\rright}{\right}
\newtheorem{theorem}{Theorem}

\newproclaim{definition}[theorem]{Definition}
\newtheorem{proposition}[theorem]{Proposition}
\newproclaim{assumption}{Assumption}
\newtheorem{lemma}[theorem]{Lemma}
\newtheorem{corollary}[theorem]{Corollary}
\newproclaim{remark}[theorem]{Remark}
\newproclaim{example}{Example}
\makeatother

\begin{document}
\begin{frontmatter}

\title{On values of repeated games with signals}
\runtitle{On values of repeated games with signals}

\begin{aug}
\author[A]{\fnms{Hugo}~\snm{Gimbert}\thanksref{m1,T1}\ead
[label=e1]{hugo.gimbert@labri.fr}},
\author[B]{\fnms{J\'er{\^o}me}~\snm{Renault}\thanksref{m2,T2}\ead
[label=e2]{jerome.renault@tse-fr.eu}},
\author[C]{\fnms{Sylvain}~\snm{Sorin}\thanksref{m3,T3}\ead
[label=e6]{sylvain.sorin@img-prg.fr}},
\author[D]{\fnms{Xavier}~\snm{Venel}\corref{}\thanksref{m4,T2,T4}\ead
[label=e4]{xavier.venel@univ-paris1.fr}}
\and
\author[E]{\fnms{Wies{\l}aw}~\snm{Zielonka}\thanksref{m5}\ead
[label=e5]{zielonka@liafa.univ-paris-diderot.fr}}
\runauthor{H. Gimbert et al.}
\affiliation{Labri\thanksmark{m1},
GREMAQ, Universit\'e Toulouse 1 Capitole\thanksmark{m2},
Sorbonne Universit\'es\thanksmark{m3},
Universit\'e Paris Diderot\thanksmark{m4}
and
Universit\'{e} Paris 1 Panth\'{e}on-Sorbonne\thanksmark{m5}}
\address[A]{H. Gimbert \\
CNRS\\
Labri\\
351 cours de la Lib\'{e}ration\\
F-33405 Talence\\
France\\
\printead{e1}}
\address[B]{J. Renault\\
TSE (GREMAQ, Universit\'e \\
\quad Toulouse 1 Capitole)\\
21 all\'ee de Brienne\\
31000 Toulouse\\
France\\
\printead{e2}\hspace*{42pt}}
\address[C]{S. Sorin\\
Sorbonne Universit\'es\\
UPMC Univ. Paris 06\\
Institut de Math\'ematiques\\
\quad de Jussieu-Paris Rive Gauche\\
UMR 7586\\
CNRS\\
Univ. Paris Diderot\\
Sorbonne Paris Cit\'e, F-75005\\
Paris\\
France\\
\printead{e6}}
\address[D]{X. Venel\\
Centre d'\'{e}conomie de la Sorbonne\\
Universit\'{e} Paris 1 Panth\'{e}on-Sorbonne\\
106-112 Boulevard de l'H\^{o}pital\\
75647 Paris Cedex 13\\
France\\
\printead{e4}}
\address[E]{W. Zielonka\\
LIAFA\\
Universit\'e Paris Diderot Paris 7\\
75205 Paris Cedex 13\\
France\\
\printead{e5}}
\end{aug}
\thankstext{T1}{Supported in part by the CNRS PEPS
project ``Jeux Stochastiques et V\'{e}rification.''}
\thankstext{T2}{Supported by the Agence Nationale
de la Recherche, under Grant ANR JEUDY, ANR-10-BLAN 0112.}
\thankstext{T3}{Supported in part by the European Union under the 7th
Framework Programme
``FP7-PEOPLE-2010-ITN'', Grant agreement number 264735-SADCO.}
\thankstext{T4}{Supported by the Israel Science Foundation under Grant
$\sharp$1517/11.}

%
\received{\smonth{7} \syear{2014}}
%
\revised{\smonth{12} \syear{2014}}

%
\begin{abstract}
We study the existence of different notions of value in two-person
zero-sum repeated games where the state evolves and players receive
signals. We provide some examples showing that the limsup value (and
the uniform value) may not exist in general. Then we show the existence
of the value for any Borel payoff function if the players observe a
public signal including the actions played.
We also prove two other positive results without assumptions on the
signaling structure: the existence of the $\sup$ value in any game and
the existence of the uniform value in recursive games with nonnegative
payoffs.
\end{abstract}

%
\end{frontmatter}

\section{Introduction}

The aim of this article is to study two-player zero-sum general
repeated games with signals (sometimes called ``stochastic
games with partial observation''). At each stage, each player chooses
some action in a finite set. This generates a stage reward then a new
state and new signals are randomly chosen 
through a transition probability depending on the current state and
actions, and with finite support. Shapley \cite{Shapley53} studied the
special case of standard stochastic games where the players observe, at
each stage, the current state and the past actions. There are several
ways to analyze these games. We will distinguish two approaches:
Borelian evaluation and uniform value.

In this article, we will mainly use a point of view coming from the
literature on determinacy of multistage games (Gale and Stewart
\cite{Gale53}). One defines a function, called \textit{evaluation}, on the
set of plays (infinite histories) and then studies the existence of a
value in the normal form game where the payoff is given by the
expectation of the evaluation, with respect to the probability induced
by the strategies of the players. Several evaluations will be considered.

In the initial model of Gale and Stewart \cite{Gale53} of two-person
zero-sum multistage game with perfect information, there is no state
variable. The players choose, one after the other, an action from a
finite set and both observe the previous choices. Given a subset $A$ of
the set of plays (in this framework: infinite sequences of actions),
player $1$ wins if and only if the actual play belongs to the set $A$:
the payoff function is the indicator function of $A$. Gale and Stewart
proved that the game is determined: either player $1$ has a winning
strategy or player $2$ has a winning strategy, in the case where $A$ is
open or closed with respect to the product topology. This result was
then extended to more and more general classes of sets until Martin
\cite{Martin75} proved the determinacy for every Borel set. When $A$
is an arbitrary subset of plays, Gale and Stewart \cite{Gale53} showed
that the game may be not determined.

In $1969$, Blackwell \cite{Blackwell69} studied the case (still
without state variable) where the players play simultaneously and are
told their choices. Due to the lag of information, the determinacy
problem is not well defined. Instead, one investigates the
probability that the play belongs to some subset $A$. When $A$ is a
$G_{\delta}$-set, a countable intersection of open sets, Blackwell
proved that there exists a real number $v$, the \textit{value} of the
game, such that for each $\varepsilon>0$, player $1$ can ensure that
the probability of the event: ``the play is in $A$'' is greater than
$v-\varepsilon$, whereas player $2$ can ensure that it is less than
$v+\varepsilon$.

The extension of this result to Shapley's model (i.e., with a state
variable) was done by Maitra and Sudderth. They focus on the specific
evaluation where the payoff is the largest stage reward obtained
infinitely often. They prove the existence of a value, called
\textit{limsup value}, in the countable framework \cite{Maitra92}, in the
Borelian framework \cite{Maitra93a} and in a finitely additive setting
\cite{Maitra93b}. In the first two cases, they assume some finiteness
of the action set (for one of the players). Their result especially
applies to finite stochastic games where the global payoff is the
limsup of the mean expected payoff.

Blackwell's existence result was generalized by Martin \cite{Martin98}
to any Borel-measurable evaluation, whereas Maitra and Sudderth \cite
{Maitra98} extended it further to stochastic games in the finitely
additive setting. In all these results, the players observe the past
actions and the current state.

Another notion used in the study of stochastic games (where a play
generates a sequence of rewards) is the uniform value where some
uniformity condition is required. Basically, one looks at the largest
amount that can be obtained by a given strategy for a family of
evaluations (corresponding to longer and longer games). There are
examples where the uniform value does not exist: Lehrer and Sorin
\cite{Lehrer92} describe such a game with a countable set of states and
only one player, having a finite action set. On the other hand,
Rosenberg, Solan and Vieille \cite{Rosenberg2002} proved the existence
of the uniform value in partial observation Markov Decision Processes
(one player) when the set of states and the set of actions are finite.
This result was extended by Renault \cite{Renault2011} to general
action space.

The case of stochastic games with standard signaling, that is, where
the players observes the state and the actions played has been treated
by Mertens and Neyman \cite{Mertens81}. They proved the existence of a
uniform value for games with a finite set of states and finite sets of
actions. In fact, their proof also shows the existence of a value for
the limsup of the mean payoff, as studied in Maitra and Sudderth and
that both values are equal.

The aim of this paper is to provide new existence results when the
players are observing only signals on state and actions. In
Section~\ref{model}, we define the model and present several specific Borel
evaluations. We then prove the existence of a value in games where the
evaluation of a play is the largest stage reward obtained along it,
called $\mathit{sup}$ evaluation and study several examples where the limsup
value does not exist.

Section~\ref{secsymmetric} is the core of this paper. We focus on the
case of symmetric signaling structure: multistage games where both
players have the same information at each stage, and prove that a value
exists for any Borel evaluation. For the proof, we introduce an
auxiliary game where the players observe the state and the actions
played and we apply the generalization of Martin's result to standard
stochastic games. Finally, in Section~\ref{uniform}, we introduce
formally the notion of uniform value and prove its existence in
recursive games with nonnegative payoffs.

\section{Repeated game with signals and Borel evaluation}\label{model}

Given a set $X$, we denote by $\Delta_f(X)$ the set of probabilities
with finite support on $X$. For any element $x\in
X$, $\delta_x$ stands for the Dirac measure concentrated on $x$.

\subsection{Model}

A \textit{repeated game form with signals} $\Gamma=(X,I,J,C,D,\pi, q)$
is defined by a set of states $X$, two finite sets of actions
$I$ and $J$, two sets of signals $C$ and $D$, an initial distribution
$\pi\in\Delta_f(X \times C \times D)$ and a transition function $q$
from $X\times I \times J$ to $\Delta_f(X\times C \times D)$. A \textit{repeated game with signals} $(\Gamma,g)$ is a pair of a repeated game
form and a reward function $g$ from $X \times I \times J$ to
$[0,1]$.

This corresponds to the general model of repeated game introduced in
Mertens, Sorin and Zamir \cite{Mertens94}.

The game is played as follows. First, a triple $(x_1,c_1,d_1)$ is drawn
according to the probability $\pi$. The initial state is $x_1$, player
$1$ learns $c_1$ whereas player~$2$ learns~$d_1$. Then, independently,
player $1$ chooses an action $i_1$ in $I$ and player~$2$ chooses an
action $j_1$ in $J$. A new triple $(x_2,c_2,d_2)$ is drawn according to
the probability distribution $q(x_1,i_1,j_1)$, the new state is $x_2$,
player $1$ learns $c_2$, player~$2$ learns $d_2 $ and so on.
At each stage $n$ players choose actions $i_n$ and $j_n$ and a triple
$(c_{n+1},d_{n+1} ,x_{n+1})$ is drawn according to
$q(x_n,i_n,j_n)$, where $x_n$ is the current state, inducing the
signals received by the players and the state at the next stage.

For each $n\geq1$, we denote by $H_n=(X \times C \times D \times I
\times J)^{n-1} \times X \times C \times D$ the set of \textit{finite
histories} of length $n$, by $H^1_n=(C \times I )^{n-1} \times C$ the
set of histories of length $n$ for player $1$ and by $H^2_n=(D \times J
)^{n-1} \times D$ the set of histories of length $n$ for player $2$.
Let\vspace*{1pt} $H = \bigcup_{n\geq1} H_n$.

Assuming perfect recall, a \textit{behavioral strategy} for player 1 is a
sequence $\sigma=(\sigma_n)_{n \geq1}$,
where $\sigma_n$, the strategy at stage $n$, is a mapping from $H^1_n$
to $\Delta(I)$, with the interpretation that $\sigma_n(h)$ is the
lottery on actions used by player 1 after $h\in H^1_n$.
In particular, the strategy $\sigma_1$ at stage $1$
is simply a mapping from $C$ to $\Delta(I)$ giving the law of the
first action played by player 1 as a function of his initial signal.
Similarly, a \textit{behavorial strategy} for player 2 is a sequence
$\tau
=(\tau_n)_{n \geq1}$, where $\tau_n$ is a mapping from $H^2_n$ to
$\Delta(J)$. We denote by $\Sigma$ and $\mathcal{T}$ the sets of behavioral
strategies of player 1 and player 2, respectively.

If for every $n\geq1$ and $h\in H^1_n$, $\sigma_n(h)$ is a Dirac
measure then
the strategy is \textit{pure}. A \textit{mixed strategy} is a
distribution over pure strategies.

Note\vspace*{1pt} that since the initial distribution $\pi$ and the transition $q$
have finite support and the sets of actions are finite, there exists a
finite subset $H^0_n \subset H_n$ such that for all strategies $(\sigma
,\tau)$ the set of histories that are reached at stage $n$ with a
positive probability is included in $H^0_n$.

Hence, no additional measurability assumptions on the strategies are
needed. It is standard that a pair of strategies $(\sigma, \tau)$
induces a probability
$\P_{\sigma, \tau} $ on the set of \textit{plays} $H_\infty=(X\times C
\times D \times I \times J)^{\infty}$ endowed with the
$\sigma$-algebra $\mathcal{H}_\infty$ generated by the cylinders
above the
elements of $H$. We denote by $\mathbb{E}_{\sigma,\tau}$ the
corresponding expectation.

Historically, the first models of repeated games assumed that both
$c_{n+1}$ and $d_{n+1}$ determine $(i_n, j_n)$ (standard signalling on
the moves also called ``full monitoring'').

A \textit{stochastic game} corresponds to the case where in addition the
state is known: both $ c_{n+1}$ and $d_{n+1}$ contain $x_{n+1}$.

A \textit{game with incomplete information} corresponds to the case where
in addition the state is fixed: $x_1= x_n, \forall n$, but not
necessarily known by the players.

Several extensions have been proposed and studied; see, for example,
Neyman and Sorin \cite{Neyman03} in particular Chapters~3, 21, 25,
28.

It has been noticed since Kohlberg and Zamir \cite{Kohlberg74b} that
games with incomplete information, when the information is \textit{symmetric}: $c_{n+1}= d_{n+1}$ and contains $(i_n, j_n)$, could be
analyzed by introducing an auxiliary stochastic game. However, the
state variable in this auxiliary stochastic game is no longer $x_n \in
X$ but the (common) conditional probability on $X$ given the signals,
that can be computed by the players: namely the law of $x_n$ in $\Delta
(X)$. Since then, this approach has been extended; see, for example,
Sorin \cite{Sorin2003b}, Ghosh et al. \cite{Gosh} and the analysis in
the current article shows that general repeated games with symmetric
information are the natural extension of standard stochastic games.

\subsection{Borel evaluation and results}

We now describe several ways to evaluate each play and the
corresponding concepts. We
follow the multistage game determinacy literature and define an
evaluation function $f$ on infinite plays. Then we study the existence
of the value of the normal form game $(\Sigma,\mathcal{T},f)$. We
will consider
especially four evaluations: the general Borel evaluation, the sup
evaluation, the limsup evaluation and the limsup-mean evaluation.

A Borel \textit{evaluation} is a $\mathcal{H}_\infty$-measurable function
from the set of plays $H_\infty$ to $[0,1]$.

\begin{definition}
Given an evaluation $f$, the game $\Gamma$ has a \textit{value} if
\[
\sup_\sigma\inf_{\tau} \mathbb{E}_{\sigma,\tau}
( f )=\inf_\tau\sup_\sigma
\mathbb{E}_{\sigma,\tau} ( f ).
\]
This real number is called the value and denoted by $v(f)$.
\end{definition}
%



Given a repeated game $(\Gamma,g)$, we will study several specific
evaluations defined through the stage payoff function $g$.

\subsubsection{Borel evaluation: $\sup$ evaluation}\label{bor2}

The first evaluation is the supremum evaluation where a play is
evaluated by the largest payoff obtained along it.

\begin{definition}
$\gamma^s$ is the \textit{sup} evaluation defined by
\[
\forall h\in H_\infty, \qquad {\gamma^s}(h)=\sup
_{n\geq1} g(x_n,i_n,j_n).
\]
In $(\Sigma,\mathcal{T},\gamma^s)$, the $\max\min$, the $\min\max
$, and the
value (called the \textit{sup value} if it exists) are, respectively,
denoted by $\underline{v}^s$,
$\overline{v}^s$ and $v^s$.
\end{definition}

The specificity of this evaluation is that for every $n\geq1$, the
maximal stage payoff obtained before $n$ is a lower bound of the
evaluation on the current play. We prove that the $\sup$ value always exists.

\begin{theorem}\label{sup}
A repeated game $(\Gamma, g)$ with the $\sup$ evaluation has a value $v^s$.
\end{theorem}

For the proof, we use the following result. We call \textit{strategic
evaluation} a function $F$ from $\Sigma\times\tau$ to $[0,1]$. It is
clear that an evaluation $f$ induces naturally a strategic evaluation
by $F(\sigma,\tau)=\mathbb{E}_{\sigma,\tau}  ( f  )$.

\begin{proposition}\label{FF}
Let $(F_n)_{n\geq1}$ be an increasing sequence of strategic
evaluations from $\Sigma\times\tau$ to $[0,1]$ that converges to some
function $F$.
Assume that:
\begin{itemize}
\item $\Sigma$ and $\tau$ are compact convex sets,
\item for every $n\geq1$, $F_n(\sigma,\cdot)$ is lower semicontinuous and
quasiconvex on $\tau$ for every $\sigma\in\Sigma$,
\item for every $n\geq1$, $F_n(\cdot,\tau)$ is upper semicontinuous and
quasiconcave on $\Sigma$ for every $\tau\in\tau$.
\end{itemize}
Then the normal form game $(\Sigma,\tau,F)$ has a value $v$.
\end{proposition}

A more general version of this proposition can be found in Mertens,
Sorin and Zamir \cite{Mertens94} (Part A, Exercise 2, Section~1.f, page~10).

\begin{pf*}{Proof of Theorem~\protect\ref{sup}}
Let $n\geq1$ and define the strategic evaluation $F_n$ by
\[
F_n(\sigma,\tau)= \mathbb{E}_{\sigma,\tau} \Bigl(\sup
_{t\leq n} g(x_t,i_t,j_t)
\Bigr).
\]
Players remember their own previous actions so by Kuhn's theorem \cite
{Kuhn53}, there is equivalence between mixed strategies and behavioral
strategies. The sets of mixed strategies are naturally convex.
The set of histories of length $n$ having positive probability is
finite and, therefore, the set of pure strategies is finite. For every
$n\geq1$, the function $F_n(\sigma,\tau)$ is thus the linear extension
of a finite game. In particular $F_n(\sigma,\cdot)$ is lower semicontinuous
and quasiconvex on $\tau$ for every $\sigma\in\Sigma$ and upper
semicontinuous and quasiconcave on $\Sigma$ for every $\tau\in\tau$.

Finally, the sequence $(F_n)_{n\geq1}$ is increasing to
\[
F(\sigma,\tau)=\mathbb{E}_{\pi,\sigma,\tau} \Bigl(\sup_{t}
g(x_t,i_t,j_t) \Bigr).
\]

It follows from Proposition~\ref{FF} that the game $\Gamma$ with the
$\sup$ evaluation has a value.
\end{pf*}

\subsubsection{Borel evaluation: $\operatorname{limsup}$ evaluation}

Several authors have especially focused on the $\operatorname{limsup}$ evaluation and
the $\operatorname{limsup}$-mean evaluation.

\begin{definition}
$\gamma^*$ is the \textit{limsup} evaluation defined by
\[
\forall h\in H_\infty, \qquad \gamma^*(h)=\limsup_n
g(x_n,i_n,j_n).
\]
In $(\Sigma,\mathcal{T},\gamma^*)$, the $\max\min$, the $\min\max
$, and the
value (called the \textit{limsup value}, if it exists) are,
respectively,
denoted by $\underline{v}^*$, $\overline{v}^*$ and $v^*$.
\end{definition}

\begin{definition}
$\gamma_m^*$ is the \textit{limsup-mean} evaluation defined by
\[
\forall h\in H_\infty,\qquad \gamma_m^*(h)=\limsup
_n \frac{1}{n}\sum_{t=1}^n
g(x_t,i_t,j_t).
\]
In $(\Sigma,\mathcal{T},\gamma_m^*)$, the $\max\min$, the $\min
\max$, and the
value (called the \textit{limsup-mean value}, if it exists) are,
respectively, denoted by $\underline{v}_m^*$, $\overline{v}_m^*$ and $v_m^*$.
\end{definition}
The limsup-mean evaluation is closely related to the limsup evaluation.
Indeed, the analysis of the limsup-mean evaluation of a
stochastic game can be reduced to the study of the limsup evaluation of
an auxiliary stochastic game having as set of states the set of finite
histories of the original game.

These evaluations were especially studied by Maitra and Sudderth \cite{Maitra92,Maitra93a}. In \cite{Maitra92}, they proved the
existence of the limsup value in a stochastic game with a countable set
of states and finite sets of actions when the players observe the state
and the actions played. Next, they extended in \cite{Maitra93a} this
result to Borel measurable evaluation.

We aim to study potential extensions of their results to repeated game
with signals. In general, a repeated game with signals has no value
with respect to the limsup evaluation as shown in the following three
examples. In each case, we also show that the limsup-mean value does
not exist.

\begin{example}\label{guess}
We consider a recursive game where the players observe neither the
state nor the action played by the other player. We say that the
players are \textit{in the dark}.

This example, due to Shmaya, is also described in Rosenberg, Solan and
Vieille \cite{Rosenberg2009} and can be interpreted as ``pick the
largest integer.''

The set of states is finite $X=\{s_1,s_2,s_3, 0^*, 1^*,-1^*,2^*,-2^*\}
$, the action set of player $1$ is $I=\{T,B\}$, the action set of
player $2$ is $J=\{L,R\}$, and the transition is given by
\[
\begin{tabular}{cccc}
& $
\begin{array}{c@{\quad}c} L & R
\end{array}
$ & $
\begin{array}{c@{\quad}c} L \hspace{14mm} & \hspace{13mm} R
\end{array}
$ & $
\begin{array}{c@{\quad}c} L & R
\end{array}
$ \\
$
\begin{array}{cc} T \\ B
\end{array}
$
&$\lleft(
\begin{array}{c@{\quad}c} s_1 & -2^* \\ s_1 & -2^*\\
\end{array}
 \rright)$ & $\lleft(
\begin{array}{c@{\quad}c} s_2 & 1/2\bigl(-1^*\bigr) +1/2 (s_3) \vspace*{2pt}\\ 1/2\bigl(1^*\bigr) +1/2 (s_1) &
0^*\\
\end{array}
 \rright)$ & $\lleft(
\begin{array}{c@{\quad}c} s_3 & s_3 \\ 2^* & 2^*\\
\end{array}
 \rright)$ \\
& $s_1$ & $ s_2$ & $ s_3 $
\end{tabular}\hspace*{-6pt}.
\]
The payoff is $0$ in states $s_1$,$s_2$, and $s_3$. For example, if the
state is $s_2$, player $1$ plays~$T$ and player $2$ plays $R$ then with
probability $1/2$ the payoff is $-1$ forever, and with probability
$1/2$ the next state is $s_3$. States denoted with a star are absorbing
states: if state $k^*$ is reached, then the state is $k^*$ for the
remaining of the game and the payoff is $k$.

\begin{cl*}\label{cl1}
The game which starts in $s_2$ has no limsup value:
$\underline{v}^*= -1/2< 1/2=\overline{v}^*$.
\end{cl*}

Since the game is recursive, the limsup-mean evaluation and the limsup
evaluation coincide, so there is no limsup-mean value either. It also
follows that the uniform value, defined formally in Section~\ref{uniform}, does not exist.

\begin{pf*}{Proof of \hyperref[cl1]{Claim}}
The situation is symmetric, so we
consider what player $1$ can guarantee.

After player $1$ plays $B$, the game is essentially over from player
$1$'s viewpoint: either absorption occurs or the state moves to $s_1$
where player $1$'s actions are irrelevant. Therefore, the only relevant
past history in order to define a strategy of player $1$ corresponds to
all his past actions being $T$. A strategy of player $1$ is thus
specified by the probability $\varepsilon_n$ to play $B$ for the first
time at stage $n$; let $\varepsilon^*$ be the probability that player
$1$ plays $T$ forever.

Player 2 can reply as follows: fix $\varepsilon>0$, and consider $N$
such that\break $\sum_{n=N}^{\infty} \varepsilon_n \leq\varepsilon$. Define
the strategy $\tau$ which plays $L$ until stage $N-1$ and $R$ at stage
$N$. For any $n>N$, we have
\[
\mathbb{E}_{s_2, \sigma, \tau} \bigl(g (x_n,i_n,j_n)
\bigr)\leq \varepsilon ^*(-1/2)+ \Biggl( \sum_{n=1}^{N-1}
\varepsilon_n \Biggr) (-1/2) + \varepsilon (1/2 ) \leq-1/2 +
\varepsilon.
\]
It follows that player 1 cannot guarantee more than $-1/2$ in the
limsup sense.
\end{pf*}
\end{example}

\begin{example}\label{semiguess}
We consider a recursive game where one player is more informed than
the other: player $2$ observes the state variable and the past actions
played whereas player $1$ observes neither the state nor the actions played.

This structure of information has been studied, for example, by
Rosenberg, Solan, and Vieille \cite{Rosenberg2004}, Renault \cite
{Renault2012a} and Gensbittel, Oliu-Barton and Venel \cite
{Gensbittel2013}. They proved the existence of the uniform value under
the additional assumption that the more informed player controls the
evolution of the beliefs of the other player on the state variable.

The set of states is finite $X=\{s_2,s_3, 0^*, 1/2^*,-1^*,2^*\}$, the
action set of player $1$ is $I=\{T,B\}$, the action set of player $2$
is $J=\{L,R\}$, and the transition is given by
\[
\begin{tabular}{c@{\quad}c@{\quad}c}
& $
\begin{array}{c@{\quad}c} L \hspace{13mm} & \hspace{7mm} R
\end{array}
$ & $
\begin{array}{c@{\quad}c} L & R
\end{array}
$\\
$
\begin{array}{cc} T \\ B
\end{array}
$ & $\lleft(
\begin{array}{c@{\quad}c}
s_2 & 1/2\bigl(-1^*\bigr) +1/2 (s_3) \vspace*{2pt}\\
(-1/2)^* & 0^*\\
\end{array}
 \rright)$ & $\lleft(
\begin{array}{c@{\quad}c}
s_3 & s_3 \\
2^* & 2^*\\
\end{array}
 \rright)$ \\
& $s_2$ & $s_3$
\end{tabular}\hspace*{-6pt}.
\]
We focus on the game which starts in $s_2$. Both players can guarantee
$0$ in the $\sup$ evaluation: player $2$ by playing $L$ forever and
player $1$ by playing $T$ at the first stage and then $B$ forever.
Since the game is recursive, the limsup-mean evaluation and the limsup
evaluation are equals.

\begin{cl*}\label{cl2}
The game which starts in $s_2$ has no
limsup value: $\underline{v}^*=-1/2<-1/6=\overline{v}^*$.
\end{cl*}

\begin{pf}
The computation of the $\max\min$ with
respect to the limsup-mean evaluation is similar to the computation of
Example~\ref{guess}. The reader can check that player $1$ cannot
guarantee more than $-1/2$.

We now prove that the $\min\max$ is equal to $-1/6$. Contrary to
Example~\ref{guess}, player $2$ observes the state and actions,
nevertheless the game is from his point of view strategically finished
as soon as $B$ or $R$ is played: if $B$ is played then absorption
occurs, if $R$ is played then either absorption occurs or the state
moves to $s_3$ where player 2's action are irrelevant. Therefore, when
defining the strategy of player $2$ at stage $n$, the only relevant
past history is $(s_2,T,L)^n$ and a strategy of player $2$ is defined
by the probability $\varepsilon_n$ that he plays $R$ for the first time
at stage $n$ and the probability $\varepsilon^*$ that he plays $L$ forever.

Fix $\varepsilon>0$, and consider $N$ such that $\sum_{n=N}^{\infty}
\varepsilon_n \leq\varepsilon$. Player 1's replies can be reduced to
the two following strategies: $\sigma_1$ which plays $T$ forever and,
$\sigma_2$ which plays $T$ until stage $N-1$ and $B$ at stage $N$. All
the other strategies are yielding a payoff smaller with an $\varepsilon
$-error. The strategy $\sigma_1$ yields $0  \varepsilon^* +
(1-\varepsilon^*)(-1/2) $ and the strategy $\sigma_2$ yields $(-1/2)
\varepsilon^* +(1-\varepsilon^*) 1/2- \varepsilon$.

The previous payoff functions are almost the payoff of the two-by-two
game where player $1$ chooses $\sigma_1$ or $\sigma_2$ and player $2$
chooses either never to play $R$ or to play $R$ at least once:
\[
\pmatrix{
0 & -1/2
\vspace*{2pt}\cr
-1/2 & 1/2}.
\]
The value of this game is $-1/6$, giving the result.
\end{pf}
\end{example}

\begin{example}\label{bigmatch}
In the previous examples, the state is not known to at least one player.

The following game is a variant of the Big Match introduced by
Blackwell and Ferguson \cite{Blackwell68}. It is an absorbing game:
every state except one are absorbing. Since there is only one state
where players can influence the transition and the payoff, the
knowledge of the state is irrelevant. Players can always consider that
the current state is the nonabsorbing state.

We assume that player $2$ observes the past actions played whereas
player $1$ does not
(in the original version, both player $1$ and player $2$ were
observing the state and past actions):
\[
\begin{array}{c@{\quad}c}
&
\hspace*{-2pt}\begin{array}{c@{\quad}c}
L & R \\
\end{array}
\\
\begin{array}{c}
T \\
B \\
\end{array}
&
\lleft(
\begin{array}{c@{\quad}c}
1^* & 0^* \\
0 & 1 \\
\end{array}
 \rright).\\
\end{array}
\]

\begin{cl*}\label{cl3}
The game with the sup evaluation has a
value $v_s=1$. The game with the limsup evaluation and the game with
the limsup-mean evaluation do not have a value: $\underline{v}^*=\underline{v}_m^*=0 < 1/2=\overline{v}_m^*=\underline{v}^*$.
\end{cl*}

\begin{pf}
We first prove the existence of the value
with respect to the sup evaluation. Player $1$ can guarantee the payoff
$1$. Let $\varepsilon>0$, and $\sigma$ be the strategy which plays $T$
with probability $\varepsilon$ and $B$ with probability $1-\varepsilon$. This
strategy yields a sup evaluation greater than $1-\varepsilon$. Since
$1$ is the maximum payoff, it is the value: $v^s=1$.

We now focus on the limsup evaluation and the limsup-mean evaluation.

After player $1$ plays $T$ absorption occurs. Therefore, the only
relevant past history in order to define a strategy of player $1$
corresponds to all his past actions being $B$. Let $\varepsilon_n$ be
the probability that player $1$ plays $T$ for the first time at stage
$n$ and $\varepsilon^*$ be the probability that player $1$ plays $B$ forever.

Player 2 can reply as follows: fix $\varepsilon>0$, and consider $N$
such that\break $\sum_{n=N}^{\infty} \varepsilon_n \leq\varepsilon$. Define
the strategy $\tau$ which plays $R$ until stage $N-1$ and $L$ at stage
$N$. For any $n>N$, we have
\[
\mathbb{E}_{s, \sigma, \tau} \bigl(g (x_n,i_n,j_n)
\bigr)\leq \varepsilon ^*0+ \Biggl( \sum_{n=1}^{N-1}
\varepsilon_n \Biggr)0 + \varepsilon (1 ) \leq\varepsilon.
\]

Let us compute what player $2$ can guarantee with respect to the limsup
evaluation. The computation is similar for the limsup-mean evaluation.
First, player $2$ can guarantee $1/2$ by playing the following mixed
strategy: with probability $1/2$, play $L$ at every stage and with
probability $1/2$, play $R$ at every stage.

We now prove that it is the best payoff that player $2$ can achieve.
Fix a strategy $\tau$ for player $2$ and consider the induced law $\P$
on the set $H_\infty=\{L,R\}^{\infty}$ of infinite sequences of $L$ and
$R$ induced by $\tau$ when player $1$ plays $B$ at every stage. Denote
by $\beta_n$ the probability that player $2$ plays $L$ at stage $n$. If
there exists a stage $N$ such that $\beta_{N} \geq1/2$, then playing
$B$ until $N-1$ and $T$ at stage $N$ yields a payoff greater than $1/2$
to player $1$. If for every $n$, $\beta_n\leq1/2$, then the stage
payoff is in expectation greater than $1/2$ when player $1$ plays $B$.
Therefore, the expected $\mathit{limsup}$ payoff is greater than $1/2$.
\end{pf}
\end{example}

\section{Symmetric repeated game with Borel evaluation} \label{secsymmetric}

Contrary to the $\sup$ evaluation, in general the existence of the
value for a given evaluation depends on the signaling structure. In
Section~\ref{model}, we analyzed three games without $\operatorname{limsup}$-mean
value. In this section, we prove that if the signaling structure is
symmetric as defined next, the value always exists in every Borel evaluation.

\subsection{Model and results}

\begin{definition}
A \textit{symmetric signaling repeated game form} is a repeated game
form with signals
$\Gamma=(X,I,J,C,D, \pi, q)$ such that there exists a set $S$ with
$C=D=I \times J \times S$ satisfying
\[
\forall(x,i,j)\in X\times I \times J, \qquad \sum_{s,x'}
q(x,i,j) \bigl(x',(i,j,s),(i,j,s)\bigr)=1
\]
and the initial distribution $\pi$ is also symmetric: $\pi(x,c,d)>0$
implies $c=d$.
\end{definition}

Intuitively, at each stage of a symmetric signaling repeated game form,
the players observe both actions played and a public signal $s$. It
will be convenient to write such a game form as a tuple $\Gamma
=(X,I,J,S, \pi, q)$ and since for such a game:
$q(x,i,j)(x',(i',j',s'),(i'',j'',s''))>0$ only if $i=i'=i''$ and
$j=j'=j''$ and $s'=s''$,
without loss of generality, we can and will write
$q(x,i,j)(x',s)$ as a shorthand for
$q(x,i,j)(x',(i,j,s),(i,j,s))$.
With this notation $q(x,i,j)$ and the initial distribution
$\pi$ are elements of $\Delta_f(X\times S)$. The set of observed plays
is then $V_\infty=(S\times I \times J)^\infty$.

\begin{theorem}\label{theo3}
Let $\Gamma$ be a symmetric signaling repeated game form. For every
Borel evaluation $f$, the game $\Gamma$ has a value.
\end{theorem}

\begin{corollary}\label{corotheo3}
A symmetric signaling repeated game $(\Gamma,g)$ has a limsup value and
a limsup-mean value.
\end{corollary}

\subsection{Proof of Theorem~\texorpdfstring{\protect\ref{theo3}}{8}}

Let us first give an outline of the proof. Given a symmetric signaling
repeated game form $\Gamma$ and a Borel evaluation $f$, we construct an
auxiliary standard stochastic game $\widehat{\Gamma}$ (where the
players observe the state and the actions) and a Borel evaluation
$\widehat{f}$ on the corresponding set of plays.\vspace*{1.5pt} We use the existence
of the value in the game $\widehat{\Gamma}$ with
respect to the evaluation $\widehat{f}$ to deduce the existence of the
value in the original game.

The difficult point is the definition of the evaluation $\widehat{f}$.
The key idea is to define a conditional probability with respect to the
$\sigma$-algebra of observed plays. For a given probability on plays,
the existence of such conditional probability is easy since the sets
involved are Polish. In our case, the difficulty comes from the
necessity to have the same conditional probability for any of the
probability distributions that could be generated by the strategies of
the players (Sections~\ref{secfinite}--\ref{sectrzy}). (As remarked
by a referee the observed plays generate in fact a sufficient statistic
for the plays with respect to all these distributions.) The definition
of the conditional probability is achieved in three steps: we first
define the conditional probability of a finite history with respect to
a finite observed history, then we use a martingale result to define
the conditional probability of a finite history with respect to an
observed play and finally we rely on Kolmogorov extension theorem to
construct a conditional probability on plays. Finally, we introduce the
function $\widehat{f}$ on the observed plays as the integral of $f$
with respect to this conditional probability.

After introducing few notations we prove the existence of the value by
defining the game $\widehat{\Gamma}$, assuming the existence of the
function $\widehat{f}$ (Lemma~\ref{transfer}). The next three sections
will be dedicated to the construction of the conditional probability,
then to the definition and properties of the function $\widehat{f}$ for any
Borelian payoff function~$f$.

Let $\Gamma$ be a symmetric signaling repeated game form, we do not
assume the Borel evaluation to be given.

\subsection{Notation}

Let $H_n=(X\times S \times I \times J)^{n-1}\times X \times S$, $H =
\bigcup_{n\geq1} H_n$, the set of histories and $H_\infty= (X\times S
\times I \times
J)^\infty$, the set of plays.

For all $h\in H_\infty$, define $ h |_{n} \in H_n$ as the
projection of $h$ on the $n$ first stages.

For all $h_n\in H_n$, denote by $h_n^+$ the cylinder generated by
$h_n$ in $H_\infty$: $h_n^+=\{h \in H_\infty, h
|_{n}=h_n \}$ and by ${\mathcal H}_n$ the corresponding
$\sigma$-algebra. ${\mathcal H}_\infty$ denotes the $\sigma$-algebra
generated by $\bigcup_n {\mathcal H}_n$.

Let $V_n=(S \times I \times J)^{n-1}\times S = H^1_n = H^2_n $, $V =
\bigcup_{n\geq1} V_n$ and $V_\infty=(S\times I \times J)^\infty$.

For all $v\in V_\infty$, define $ v |_{n} \in V_n$ as the
projection of $v$ on the $n$ first stages.

For all $v_n\in V_n$, denote by $v_n^+$ the cylinder generated by
$v_n$ in $V_\infty$: $v_n^+=\{v \in V_\infty, v
|_{n}=v_n \}$ and by ${\mathcal V}_n$ the corresponding
$\sigma$-algebra. ${\mathcal V}_\infty$ is the $\sigma$-algebra
generated by $\bigcup_n {\mathcal V}_n$.

We denote by $\Theta$ the application from $H_\infty$ to $V_\infty$
which forgets all the states: more precisely,
$\Theta( x_1, s_1, i_1, j_1,\ldots, x_n, s_n, i_n, j_n, \ldots) = (
s_1, i_1, j_1, \ldots,\break  s_n, i_n, j_n,\ldots)$.
We use the same notation for the corresponding application
defined from $H$ to $V$.

We denote by ${\mathcal V}^*_n$ (resp., ${\mathcal V}^*_\infty$) the
image of ${\mathcal V}_n$ (resp., ${\mathcal V}_\infty$)
by $\Theta^{-1}$ which are sub $\sigma$-algebras of ${\mathcal H}_n$
(resp., ${\mathcal H}_\infty$).
Explicitly, for $v_n\in V_n$, $v_n^*$ denotes the cylinder generated by
$v_n$ in $H_\infty$: $v_n^*=\{h\in H_\infty, \Theta(h)
|_{n}=v_n \}$, ${\mathcal V}^*_n$ are the corresponding
$\sigma$-algebras and ${\mathcal V}^*_\infty$ the $\sigma$-algebra
generated by their union.

Any ${\mathcal V}_n$ (resp., ${\mathcal V}_\infty$)-measurable function
$\ell$ on $V_\infty$ induces a ${\mathcal V}^*_n$ (resp., ${\mathcal
V}^*_\infty$)-measurable function $\ell\circ\Theta$ on $H_\infty$.

Define $\alpha$ from $H$ to $[0,1]$ where for $h_n = ( x_1, s_1, i_1,
j_1,\ldots, x_n, s_n)$:
\[
\alpha( h_n) =\pi(x_1,s_1)
\prod_{t=1}^{n-1} q(x_t,i_t,j_t)
(x_{t+1},s_{t+1})
\]
and $\beta$ from $V$ to $[0,1]$ where for $v_n = (s_1,
i_1, j_1,\ldots, s_n)$:
\[
\beta(v_n)= \sum_{ h_n \in H_n ; \Theta(h_n) = v_n }
\alpha(h_n).
\]

Let
${\overline H}_n = \{ h_n \in H_n$; $ \alpha( h_{n}) >0 \}$ and
${\overline V}_n = \Theta({\overline H}_n)$ and recall that these sets
are finite.
We introduce now the set of plays and observed plays that can occur
during the game as ${\overline H}_\infty=
\bigcap_n \overline{H}_n ^+$ and ${\overline V}_\infty= \Theta
({\overline H}_\infty) = \bigcap_n \overline{V}_n$. Remark that both are
measurable subsets of $H_\infty$ and $V_\infty$, respectively.

For every pair of strategies $(\sigma,\tau)$, we denote by $\P
_{\sigma
,\tau}$ the probability distribution induced over
the set of plays $(H_\infty, {\mathcal H}_\infty)$ and by $\mathbb
{Q}_{\sigma
,\tau}$ the probability distribution over the set of
observed plays $(V_\infty, {\mathcal V}_\infty)$. Thus, $\mathbb
{Q}_{\sigma,\tau
}$ is the image of $\P_{\sigma,\tau}$ under~$\Theta$. Note that $\operatorname{supp}( \P_{\sigma,\tau} )
\subset{\overline
H}_\infty$. We denote, respectively, by $\mathbb{E}_{\P_{\sigma,\tau
}}$ and $\mathbb{E}
_{\mathbb{Q}_{\sigma,\tau}}$ the corresponding expectations.

It turns out that for technical reasons
it is much more convenient to work with the space
$\overline{V}_\infty$ rather than with $V_\infty$ (and with
$\overline
{H}_\infty$ rather than with $H_\infty$).
And then, abusing slightly the notation,
$\mathcal{V}_\infty$ and $\mathcal{V}_n$ will tacitly denote
the restrictions to $\overline{V}_\infty$
of the corresponding $\sigma$-algebras defined on $V_\infty$.
On rare occasions this can lead to a confusion and then
we will write, for example, $\overline{\mathcal{V}}_n$ to denote the
$\sigma$-algebra $\{ U\cap\overline{V}_\infty\vert  U\in
\mathcal{V}_n \}$ the restriction of $\mathcal{V}_n$ to $\overline
{V}_\infty$.

\subsubsection{Definition of an equivalent game}\label{secconclusion}

Let\vspace*{1pt} us define an auxiliary stochastic game $\widehat{\Gamma}$. The sets
of actions $I$ and $J$ are the same as in $\Gamma$. The set of states is
$V=\bigcup_{n\geq1} V_n$ and the transition $\widehat{q}$ from
$V\times I
\times J$ to $\Delta(V)$ is given by
\[
\forall v_n\in V_n, \forall i\in I, \forall j\in J,
\qquad \widehat{q}(v_n,i,j)=\sum_{s \in S}
\psi(v_n, i, j, s)\delta_{v_n, i,
j, s},
\]
where\vspace*{1pt} $\psi(v_n, i, j, s)= \frac{\beta(v_n, i, j,
s)}{\beta(v_n)}$.

Note that if $v_n\in V_n$ then the support of $\widehat{q}(v_n,i,j)$ is
included in $V_{n+1}$, in particular is finite.
Moreover, if $\widehat{q}(v_n,i,j)(v_{n+1})>0$ then $v_{n+1}|_n=v_n$.
The initial distribution of $\widehat{\Gamma}$ is the marginal
distribution $\pi^S$ of $\pi$ on $S$, if $s\in S= V_1$, then $\pi
^S(s)=\sum_{x\in X}\pi(x,s)$ and $\pi^S(v)=0$ for $v\in V\setminus V_1$.

Let us note that the original game $\Gamma$ and the auxiliary game
$\widehat{\Gamma}$ have the same sets of strategies. Indeed a
behavioral strategy in $\Gamma$ is a mapping from $V$ to probability
distributions over actions. Thus, each behavioral strategy in $\Gamma$
is a stationary strategy in $\widehat{\Gamma}$. On the other hand
however, each state of $\widehat{\Gamma}$ ``contains''
all previously visited states and all played actions; thus, for all
useful purposes,
in $\widehat{\Gamma}$ behavioral strategies and stationary strategies
coincide.

Now suppose that $(v_1, i_1, j_1, v_2,i_2, j_2, \ldots)$ is a play in
$\widehat{\Gamma}$. Then $v_{n+1}|_n=v_n$ for all $n$ and there exists
$v\in V_\infty$ such that $v|_n=v_n$ for all $n$.
Thus, defining a payoff on infinite histories in $\widehat{\Gamma}$
amounts to defining a payoff on $V_\infty$.

\begin{lemma}\label{transfer}
Given a Borel function $f$ on $H_\infty$, there exists a Borel function
$\widehat{f}$ on $V_\infty$ such that
%
\begin{equation}
\label{eqfinal}
\mathbb{E}_{\P_{\sigma,\tau}}(f)=\mathbb{E}_{\mathbb{Q}_{\sigma
,\tau}} (
\widehat{f} ).
\end{equation}
\end{lemma}

Therefore, playing in $\Gamma$ with strategies
$(\sigma,\tau)$ and payoff $f$ is the same as playing in $\widehat
{\Gamma}$ with
the same strategies and payoff $\widehat{f}$.

By Martin~\cite{Martin98} or Maitra and Sudderth~\cite{maitra2003},
the stochastic game
$\widehat{\Gamma}$ with payoff $\widehat{f}$ has a value implying that
$\Gamma$ with payoff $f$ has the same value, which completes the proof
of Theorem~\ref{theo3}.

The three next sections are dedicated to the proof of Lemma~\ref{transfer}.

\subsubsection{Regular conditional probability of finite time
events with respect to finite observed histories}\label{secfinite}

For $m\geq n \geq1$, we define $\Phi_{n,m}$ from $ H_{\infty} \times
{\overline V}_\infty$ to $[0, 1 ] $ by
\[
\Phi_{n,m} ( h, v) =
\cases{ \displaystyle
\frac{ \sum_{h', h'|_n = h|_n, \Theta(h'|_m )= v
|_m}\alpha( h' |_{m})}{ \beta( v |_{m}) },& \quad\mbox{if }$\Theta( h |_{n})= v
|_{n}$,\vspace*{3pt}
\cr
0,& \quad\mbox{otherwise}. }
\]
This corresponds to the joint probability of the players on the
realization of the history $h$ up to stage $n$, given the observed
history $v$ up to stage $m$.

Since\vspace*{1pt} $\Phi_{n,m} ( h, v)$ depends only on $h|_n$ and $v|_m$, we can
see $\Phi_{n,m}$ as a function defined
on $H_n \times\overline{V}_m $ and note that its support is included
in $\overline{H}_n \times\overline{V}_m$.
On the other hand, since each set $U\in\mathcal{H}_n$ is a finite union
of cylinders $h_n^+$ for $h_n\in H_n$ such that $h_n^+\subset U$, $\Phi
_{n,m}$ can be seen as
a mapping from $\mathcal{H}_n\times\overline{V}_\infty$ into $[0,1]$,
where $\Phi_{n,m}(U,v)=\sum_{h_n, h_n^+\subseteq U}\Phi_{n,m}(h_n,v)$.
Bearing this last observation in mind, we have the following.

\begin{lemma}\label{lemkernel}
For every $m\geq n\geq1$, $\Phi_{n,m}$ is a probability kernel from
$(\overline{V}_\infty, {\mathcal V}_m)$ to $(H_\infty, {\mathcal H}_n)$.
\end{lemma}

\begin{pf}
Since $\sum_{h_n\in H_n}\Phi_{n,m}(h_n,v)=1$ for $v\in\overline
{V}_\infty$, $ \Phi_{n,m} ( \cdot,v)$ defines a probability on
$\mathcal{H}_n$.
Moreover, for any $U\in\mathcal{H}_n$, $\Phi_{n,m}( U, v)$ is a
function of the $m$ first components of $v$ hence is ${\mathcal V}_m$-measurable.
\end{pf}

\begin{lemma}\label{link}
Let $m\geq n \geq1$ and $(\sigma,\tau)$ be a pair of strategies.
Then, for every $v_m\in\overline{V}_m$ such that $\mathbb{Q}_{\sigma
,\tau
}(v_m^+)=\P_{\sigma,\tau}(v_m^*)>0$,
and every $h_n\in H_n$:
\[
\P_{\sigma,\tau}\bigl(h_n^+|v_m^*\bigr)=
\Phi_{n,m}(h_n,v_m).
\]
\end{lemma}

\begin{pf}
Let $v_m = ( s_1, i_1, j_1, \ldots, s_m)$ and $h_n \in H_n$,
\begin{eqnarray*}
&&\hspace*{-3pt} \P_{\sigma,\tau}\bigl(h_n^+|v_m^*\bigr)\\
&&\hspace*{-3pt}\qquad =
\frac{\P_{\sigma,\tau}(h_n^+ \cap v_m^*)}{ \P_{\sigma,\tau
}(v_m^*) }
\\
&&\hspace*{-3pt}\qquad=
\cases{ \displaystyle\frac{ \sum_{h', h'|_n = h_n, \theta( h'|_m) = v_m }%
{\alpha( h' |_m ) }%
W(i_1,j_1,\ldots,j_{m-1})}%
{\beta(v_m) W(i_1,j_1,\ldots,j_{m-1})}, & \hspace*{-4pt}\quad$\mbox{if }
\Theta(h _{n})= v_m |_{n}$,\vspace*{3pt}
\cr
0,&
\hspace*{-4pt}$\quad\mbox{otherwise}$,}
\end{eqnarray*}
where $W(i_1,j_1,\ldots,j_{m-1})=\prod_{t\leq
m-1}\sigma(v_m|_t)(i_t)\tau(v_m|_t)(j_t)$.
After simplification, we recognize on the right the definition of
$\Phi_{n,m}(v_m,h_n)$.
\end{pf}

We deduce the following lemma.

\begin{lemma}\label{phi-m-n}
For every pair of strategies $(\sigma,\tau)$, each $W\in\overline
{\mathcal{V}}_m$ and $U\in\mathcal{H}_n$ we have
%
\begin{equation}
\label{eqmn} \P_{\sigma,\tau}\bigl(U \cap\Theta^{-1}(W)\bigr)=\int
_W \Phi_{n,m}(U, v)\mathbb{Q}_{\sigma,\tau}(dv) .
\end{equation}
\end{lemma}
\begin{pf}
Clearly, it suffices to prove (\ref{eqmn}) for cylinders $U=h_n^+$ and
$W=v_m^+$
with $\beta(v_m)>0$.

We have
\begin{eqnarray*}
\int_{v_m^+} \Phi_{n,m}(h_n,v)
\mathbb{Q}_{\sigma,\tau}(dv) &=& \Phi_{n,m}(h_n,v_m)
\mathbb{Q}_{\sigma,\tau}\bigl(v_m^+\bigr)
\\
&=& \P_{\sigma,\tau}\bigl(h_n^+ | v_m^*\bigr)
\mathbb{Q}_{\sigma,\tau}\bigl(v_m^+\bigr)
\\
&=& \P_{\sigma,\tau}\bigl(h_n^+ | v_m^*\bigr)
\P_{\sigma,\tau}\bigl(v_m^*\bigr)
\\
&=& \P_{\sigma,\tau}\bigl(h_n^+ \cap v_m^*\bigr).
\end{eqnarray*}
\upqed
\end{pf}
Note that (\ref{eqmn}) can be equivalently written as: for every pair
of strategies $(\sigma,\tau)$,
each $W^*\in\overline{\mathcal{V}}^*_m$ and $U\in\mathcal{H}_n$
%
\begin{equation}
\label{eqmn^*} \P_{\sigma,\tau}\bigl(U \cap W^*\bigr)=\int_{W^*}
\Phi_{n,m}\bigl(U, \Theta(h)\bigr)\P_{\sigma,\tau}(dh).
\end{equation}

\subsubsection{Regular conditional probability of finite time events
with respect
to infinite observed histories} \label{secdwa}

In this paragraph, we prove that instead of defining one application
$\Phi_{n,m}$ for every pair $(m,n) $ such that $m\geq n\geq1$, one can
define a unique probability kernel $\Phi_n$ from $(\Omega_n,
{\mathcal
V}_\infty)$ to $(H_\infty,{\mathcal H}_n)$, with $\mathbb{Q}_{\sigma
,\tau
}(\Omega_n)=1$, for all $(\sigma, \tau)$, such that the extension of
Lemma~\ref{phi-m-n} holds.

For $h\in H_\infty$, let
\[
\Omega_{h} = \bigl\{ v\in\overline{V}_\infty\vert\mbox{$\Phi_{n,m}(h,v)$ converges as $m\uparrow\infty$} \bigr\}.
\]
The domain $\Omega_{h}$ is measurable (see Kallenberg~\cite
{Kallenberg97},  page~6, e.g.). Recall that $\Omega_{h}$ depends only on
$h|_n$ and write also $\Omega_{h|_n}$ for $ \Omega_{h}$. Let then
\[
\Omega_n = \bigcap_{h_n\in H_n}
\Omega_{h_n}.
\]
We define $\Phi_n \dvtx  H_\infty\times\overline{V}_\infty\to[0,1]$ by
$\Phi_n = \lim_{m\rightarrow\infty} \Phi_{n,m}$ on $ H_\infty
\times
\Omega_n $ and $0$ otherwise. As a limit of a sequence of measurable
mappings $\Phi_n$ is measurable (see Kallenberg~\cite{Kallenberg97},
page~6, e.g.).

\begin{lemma}\label{phi-n}
\textup{(i)} For each pair of strategies $(\sigma,\tau)$,
$\mathbb{Q}_{\sigma,\tau}(\Omega_n)=1$.

\textup{(ii)}
For each $v\in\Omega_n$, $\sum_{h_n\in H_n} \Phi_n(h_n, v)=1$.

\textup{(iii)}
For each $U\in\mathcal{H}_n$
the mapping $v \mapsto\Phi_n(U, v)$ is
a measurable mapping from $(\overline{V}_\infty,\mathcal{V}_\infty
)$ to
$\mathbb{R}$.

\textup{(iv)}
For each pair of strategies $(\sigma,\tau)$, for each $U\in\mathcal{H}_n$
and each $W\in\mathcal{V}_\infty$
%
\begin{equation}
\label{eqn} \P_{\sigma,\tau}\bigl(U\cap\Theta^{-1}(W)\bigr)=\int
_W \Phi_n(U,v) \mathbb{Q}_{\sigma,\tau}(dv).
\end{equation}
\end{lemma}

\begin{pf}
(i)
For $h_n\in H_n$ and
each pair of strategies $\sigma,\tau$
we define on $H_\infty$ a sequence of random variables $Z_{h_n,m}$,
$m\geq
n$,
\[
Z_{h_n,m} = \P_{\sigma,\tau} \bigl[ h_n^+ |
\mathcal{V}_m^* \bigr].
\]

As a conditional expectation of a bounded random variable with respect
to an increasing sequence of $\sigma$-algebras, $Z_{h_n,m}$ is a
martingale (with respect to
$\P_{\sigma,\tau}$), hence converges $\P_{\sigma,\tau}$-almost surely
and in $L^1$ to the random variable $Z_{h_n}=\P_{\sigma,\tau}[ h_n^+ |
\mathcal{V}_\infty^*]$.

For $m\geq n$, we define the mappings $\psi_{n,m}[h_n] \dvtx  \overline
{H}_\infty\to[0,1]$,
\[
\psi_{n,m}[h_n](h)= \Phi_{n,m}
\bigl(h_n,\Theta(h)\bigr).
\]

Let us show that for each $h_n\in H_n$, $\psi_{m,n}[h_n]$ is a version
of the conditional expectation $\mathbb{E}_{\P_{\sigma,\tau
}}[\mathbh{1}_{h_n}| \mathcal
{V}_m^*] = \P_{\sigma,\tau}  [ h_n^+ | \mathcal
{V}_m^* ]$.
First note that $\psi_{n,m}[h_n]$ is $(H_\infty,\mathcal{V}_m^*)$
measurable. Lemma~\ref{link} implies that, for $h\in\operatorname
{supp}( \P_{\sigma
,\tau} ) \subset\overline{H}_\infty$,
$\psi_{n,m}[h_n](h)=\Phi_{n,m}(h_n,\Theta(h))=\P_{\sigma,\tau}(h_n^+
|
v|_m^*)=\P_{\sigma,\tau}(h_n^+ | \mathcal{V}_m^* )(h)$, where
$v=\Theta(h)$. Hence, the claim.

Since $\psi_{n,m}[h_n]$ is a version of $\P_{\sigma,\tau}(h_n^+ |
\mathcal{V}_m^* )$, its limit
$\psi_n[h_n]$ exists and is a version of $\P_{\sigma,\tau}(h_n^+ |
\mathcal{V}_\infty^* )$, $\P_{\sigma,\tau}$-almost surely. In particular,
\begin{enumerate}[(C1)]
\item[(C1)]
the set $\Theta^{-1}(\Omega_{h_n}) = \{ h\in H_\infty\vert \mbox
{$\lim_m \psi_{n,m}[h_n](h)$ exists} \}$
is $\mathcal{V}_\infty^*$ measurable and has $\P_{\sigma,\tau
}$-measure $1$,
\item[(C2)]
for each $W^*\in\mathcal{V}_\infty^*$, $\int_{W^*} \psi_n[h_n](h)
\P_{\sigma,\tau}(dh)=
\int_{W^*} \mathbb{E}[\mathbh{1}_{h_n^+} | \mathcal{V}_\infty^*]
\P_{\sigma,\tau}=\break 
\P_{\sigma,\tau}(W^*\cap h_n^+)$.
\end{enumerate}
Note that (C1) implies that $\mathbb{Q}_{\sigma,\tau}(\Omega
_n)=1$.

(ii)
If $v\in\Omega_n$ then, for all $h_n\in H_n$, $\Phi_{n,m}(h_n,v)$
converges to
$\Phi_n(h_n, v)$. But, by Lemma~\ref{lemkernel}, $\sum_{h_n\in H_n}
\Phi_{n,m}(h_n,v)=1$.
The\vspace*{1pt} sum being with finitely many nonzero terms one\vspace*{1pt} has $\sum_{h_n\in
H_n} \Phi_n(h_n, v)=1$.

(iii) Was proved before the lemma.

(iv)
Since
$\int_W \Phi_n(h_n,v) \mathbb{Q}_{\sigma,\tau}(dv) =
\int_{\Theta^{-1}(W)} \psi_n[h_n](h) \P_{\sigma,\tau}(dh)$ for $W
\in
\mathcal{V}_\infty$, using
(C2) we get
\[
\P_{\sigma,\tau}\bigl(h_n^+\cap\Theta^{-1}(W)\bigr)=\int
_W \Phi_n(h_n,v)
\mathbb{Q}_{\sigma,\tau}(dv)
\]
for $U\in\mathcal{V}_\infty$.
\end{pf}

\subsubsection{Regular conditional probability of infinite time events
with respect
to infinite observed histories} \label{sectrzy}

In this section, using Kolmogorov extension theorem we construct from
the sequence $\Phi_n$ of probability kernels from $(\Omega_n,
{\mathcal V}_\infty)$ to $(H_\infty,{\mathcal H}_n)$, one
probability kernel $\Phi$ from $(\Omega_\infty, {\mathcal V}_\infty
)$ to $(H_\infty,{\mathcal H}_n)$, with $\mathbb{Q}_{\sigma,\tau
}(\Omega_\infty
) =1$, for all $(\sigma, \tau)$.

\begin{lemma}\label{phi}
There exists a measurable subset $\Omega_\infty$ of $V_\infty$ such
that, for all strategies $\sigma,\tau$:
\begin{itemize}
\item
$\mathbb{Q}_{\sigma,\tau}(\Omega_\infty)=1$ and
\item
there exists a probability kernel $\Phi$ from $(\Omega_\infty,
{\mathcal V}_\infty)$ to $(H_\infty,{\mathcal H}_\infty)$ such that
for each $W\in\mathcal{V}_\infty$ and $U\in\mathcal{H}_\infty$
%
\begin{equation}
\label{eqfinall} \P_{\sigma,\tau}\bigl(U\cap\Theta^{-1}(W)\bigr)=\int
_W \Phi(U,v) \mathbb{Q}_{\sigma,\tau}(dv).
\end{equation}
\end{itemize}
\end{lemma}

Before proceeding to the proof, some remarks are in order.

A probability kernel having the property given above is called a
regular conditional probability.

For given strategies $\sigma$ and $\tau$, the existence of a transition
kernel $\kappa_{\alpha,\beta}$ from $(V_\infty, {\mathcal V}_\infty)$
to $(H_\infty,{\mathcal H}_\infty)$ such that for each $U\in\mathcal
{V}_\infty$ and $A\in\mathcal{H}_\infty$
\[
\P_{\sigma,\tau}\bigl(A\cap\Theta^{-1}(U)\bigr)=\int
_U \kappa_{\sigma,\tau}(A,v) \mathbb{Q}_{\sigma,\tau}(dv)
\]
is well known provided that $V_\infty$ is a Polish space and $\mathcal
{V}_\infty$ is the Borel $\sigma$-algebra. In the current framework it
is easy to introduce an appropriate metric on $V_\infty$ such that
this condition is satisfied thus the existence of $\kappa_{\sigma
,\tau
}$ is immediately assured.

The difficulty in our case comes from the fact that we look for a
regular conditional probability
which is \textit{common for all probabilities} $\P_{\sigma,\tau}$, where
$(\sigma,\tau)$ range over all strategies of both players.

\begin{pf*}{Proof of Lemma~\protect\ref{phi}}
We follow the notation of the proof of Lemma~\ref{phi-n} and define
$\Omega_\infty=\bigcap_{n\geq1} \Omega_n$. Let $(\sigma,\tau)$
be a
couple of strategies. For every \mbox{$n\geq1$}, $\mathbb{Q}_{\sigma,\tau
}(\Omega
_n)=1$, hence $\mathbb{Q}_{\sigma,\tau}(\Omega_\infty)=1$. By
Lemma~\ref{phi-n}(ii), given $v\in\Omega_\infty$, the sequence $\{\Phi
_n(\cdot,
v) \}_{n\geq1}$
of probabilities on $\{(H_\infty, {\mathcal H}_n)\}_{n\geq1}$ is well
defined. Let us show that this
sequence satisfies the condition of Kolmogorov's extension theorem.

In fact $\Phi_{n,m}(\cdot, v)$ is defined on the power set of $H_{n}$ by
\[
\forall A \subset H_{n}, \qquad \Phi_{n,m}(A, v)=\sum
_{h_n\in A} \Phi_{n,m}(h_n, v).
\]
Thus, for every $h_n \in H_n$, we have
\begin{eqnarray*}
\Phi_{n,m}(h_n, v)&=&\frac{\P_{\sigma,\tau}(v|_m^* \cap
h_n^+)}{\P_{\sigma,\tau}(v|_m^*)}
\\
&=& \frac{\P_{\sigma,\tau}(v|_m^* \cap(h_n\times I \times J \times
X\times S)^+)}{\P_{\sigma,\tau}(v|_m^*)}
\\
&=&\Phi_{n+1,m}\bigl(h_n \times(I \times J \times X\times S),
v\bigr).
\end{eqnarray*}
Taking the limit, we obtain the same equality for $\Phi_n$ and $\Phi
_{n+1}$ hence the compatibility condition.
By the Kolmogorov extension theorem for each $v\in\Omega$,
there exists a measure $\Phi(\cdot, v)$
on $(H_\infty,\mathcal{H}_\infty)$ such that
\[
\Phi\bigl(h_n^+, v\bigr)= \Phi_n\bigl(h_n^+,
v\bigr)
\]
for each $n$ and each $h_n\in H_n$.

Let us prove that, for each $U\in\mathcal{H}_\infty$, the mapping
$v\mapsto\Phi(U,v)$ is $\mathcal{V}_\infty$-measurable on $\Omega
_\infty$.

Let $\mathcal{C}$ be the class of sets $A \in\mathcal{H}_\infty$ such
that $\Phi(A,\cdot)$ has this property. By Lemma~\ref{phi-n},
$\mathcal
{C}$ contains the $\pi$-system
consisting of cylinders generating $\mathcal{H}_\infty$. To show that
$\mathcal{H}_\infty\subseteq\mathcal{C}$ it suffices to show that
$\mathcal{C}$ is a $\lambda$-system. Let $A_i$ be an increasing
sequence of sets belonging to
$\mathcal{C}$. Since, for each $v\in\overline{V}_\infty$, $\Phi
(\cdot
,v)$ is a measure, we have $\Phi(\bigcup_n A_n,v)=\sup_n \Phi(A_n,v)$.
However, $v\mapsto\sup_n \Phi(A_n,v)$ is measurable as a supremum of
measurable mappings $v\mapsto\Phi(A_n,v)$.
Let $A\supset B$ be two sets belonging to $\mathcal{C}$. Then $\Phi
(A\setminus B,v) + \Phi(B,v)=\Phi(A,v)$ by additivity of measure and
$v\mapsto\Phi(A\setminus B,v)= \Phi(A,v) - \Phi(B,v)$ is measurable as
a difference of measurable mappings.

To prove (\ref{eqfinall}), take a measurable subset $W$ of $\overline
{V}_\infty$ and consider the set function
\[
\mathcal{H}_\infty\ni U \mapsto\int_W \Phi(U,dv)
Q_{\sigma,\tau}(dv).
\]
Since $\Phi(\cdot,v)$ is nonnegative this set function is a measure on
$(H_\infty,\mathcal{H}_\infty)$. However, by Lemma~(\ref{phi-n}), this
mapping is equal to $U \mapsto\P_{\sigma,\tau}( U \cap\Theta^{-1}(W))$
for $U$ belonging to the $\pi$-system of cylinders generating
$\mathcal
{H}_\infty$. But two measures equal on a generating $\pi$-system are
equal, which terminates the proof of (\ref{eqfinall}).
\end{pf*}

A standard property of probability kernels and the fact that $\Omega
_\infty$ has measure $1$ imply:

\begin{corollary}\label{cortransfer}
Let $f : H_\infty\to[0,1]$ be $\mathcal{H}_\infty$-measurable mapping.
Then the mapping $\widehat{f} : V_\infty\to[0,1]$ defined by
\[
\widehat{f}(v) = %
\cases{ \displaystyle\int_{H_\infty}
f(h)\Phi(dh, v), & \quad\mbox{if $v\in\Omega_\infty$},\vspace*{3pt}
\cr
0,
& \quad$\mbox{otherwise}$,}
\]
is $\mathcal{V}_\infty$-measurable and
\[
\mathbb{E}_{\P_{\sigma,\tau}}[ f ] = \mathbb{E}_{\mathbb
{Q}_{\sigma,\tau}}[ \widehat{f} ]\qquad
\forall\sigma, \tau.
\]
\end{corollary}

\begin{remark}
In the previous proof, we proceeded through a reduction from a
symmetric repeated game to a stochastic game in order to apply Martin's
existence result. The same procedure can be applied for $N$-player
repeated games. Let us consider a $N$-player symmetric signaling
repeated game. One defines a conditional probability and therefore
associates to all Borel payoffs $f^i$ on plays, $i \in N$ an associated
Borel evaluation $\widehat{f}^i$ on the space of observed plays, therefore,
reducing the problem to a $N$-player stochastic game with Borelian payoffs.

For example, Mertens  \cite{Mertens86} showed the existence
of pure $\varepsilon$-Nash equilibrium in $N$-person stochastic games
with Borel payoff functions where at each stage at most one of the
players is playing. Using the previous reduction, one can deduce the
existence of
pure $\varepsilon$-Nash equilibrium in $N$-person symmetric repeated
games with Borel payoff functions where at each stage at most one of
the players is playing.
\end{remark}

\section{Uniform value in recursive games with nonnegative
payoffs}\label{uniform}

In Section~\ref{model} and Section~\ref{secsymmetric}, we focused on
Borel evaluations. In this last section, we focus on the family of mean
average of the $n$ first stage rewards and the corresponding uniform value.

\begin{definition}
For each $n \geq1$, the \textit{mean expected payoff} induced by
$(\sigma,\tau)$ during the first $n$ stages is
\[
\gamma_n (\sigma, \tau)=\mathbb{E}_{ \sigma, \tau} \Biggl(
\frac
{1}{n}\sum_{t=1}^n g
(x_t,i_t,j_t) \Biggr).
\]
\end{definition}

\begin{definition}\label{stocuni1} \label{stocuni}
Let $v$ be a real number.

A strategy $\sigma^*$ of player $1$ \textit{guarantees} $v$ \textit{in the
uniform sense} in $(\Gamma, g)$ if for all $\eta>0$ there exists $n_0
\geq1$ such that
%
\begin{equation}
\forall n\geq n_0, \forall\tau\in\mathcal{T},\qquad \gamma
_n\bigl(\sigma^*,\tau\bigr) \geq v-\eta.
\end{equation}
Player $1$ can \textit{guarantee} $v$ \textit{in the uniform sense} in
$(\Gamma, g)$ if for all $\varepsilon>0$ there exists a strategy
$\sigma^*\in\Sigma$ which guarantees $v-\varepsilon$ in the uniform
sense.

A symmetric notion holds for player $2$.
\end{definition}

\begin{definition}
The \textit{uniform $\max\min$}, denoted by $\underline{v}_\infty$, is
the supremum of all the payoff that player $1$ can guarantee in the
uniform sense. A \textit{uniform $\min\max$} denoted by $\overline
{v}_\infty$ is defined in a dual way.

If both players can guarantee $v$ in the uniform sense, then $v$ is the
\textit{uniform value} of the game $(\Gamma, g)$ and denoted by
$v_\infty$.
\end{definition}

Many existence results have been proven in the literature concerning
the uniform value and uniform $\max\min$ and $\min\max$; see, for
example, Mertens, Sorin and Zamir \cite{Mertens94} or Sorin \cite
{Sorin2002}. Mertens and Neyman \cite{Mertens81} proved that in a
stochastic game with a finite state space and finite actions spaces,
where the players observe past payoffs and the state, the uniform value
exists. Moreover, the uniform value is equal to the limsup-mean value
and for every $\varepsilon>0$ there exists a strategy which guarantees
$v_\infty-\varepsilon$ both in the limsup-mean sense and in the
uniform sense.

In general, the uniform value does not exist (either in games with
incomplete information on both sides or in stochastic games with
signals on the actions) and in particular its existence depends upon
the signaling structure.

\begin{remark}
For $n\geq1$, the \textit{$n$-stage game} $(\Gamma_n, g)$ is
the zero-sum game with normal form $(\Sigma, \mathcal{T}, \gamma_n)$
and value $v_n$.
It is interesting to note that in the special case of symmetric
signaling repeated games with a finite set of states and finite set of
signals, a uniform value may not exist, since even the sequence of
values $v_n $ may not converge (Ziliotto \cite{Ziliotto2013}), but
there exists a value for any Borel evaluation by Theorem~\ref{theo3}.
\end{remark}

We focus now on the specific case of recursive games with nonnegative
payoff defined as follows.

\begin{definition}
Recall that a state is \textit{absorbing} if the probability to stay
in this state is 1 for all actions and the payoff is also independent
of the actions played. A repeated game is \textit{recursive} if the
payoff is equal to $0$ outside the
absorbing states. If all absorbing payoffs are nonnegative, the
game is \textit{recursive} and \textit{nonnegative}.
\end{definition}

Solan and Vieille \cite{Solan2002} have shown the existence of a
uniform value in nonnegative recursive games where the players observe
the state and past actions played. We show that the result is true
without assumption on the signals to the players.

In a recursive game, the limsup-mean evaluation and the limsup
evaluation coincide.
If the recursive game has nonnegative payoffs, the sup evaluation, the
limsup evaluation and the limsup-mean evaluation both coincide. So,
Theorem~\ref{sup} implies the existence of the value with respect to
these evaluations. Using a similar proof, we obtain the stronger theorem.

\begin{theorem}\label{recursive} A recursive game with nonnegative
payoffs has a uniform value~$v_\infty$, equal to the sup value and the
limsup value. Moreover, there exists a strategy of player $2$ that
guarantees $v_\infty$.
\end{theorem}

The proof of the existence of the uniform value is similar to the proof
of Proposition~\ref{FF} while using a specific sequence of strategic
evaluations.

\begin{pf*}{Proof of Theorem~\protect\ref{recursive}}
The sequence of stage payoffs is nondecreasing on each history:
$0$ until absorption occurs and then constant, equal
to some nonnegative real number. In particular, the payoff converges
and the $\operatorname{limsup}$ can be replaced by a limit.

Let $\sigma$ be a strategy of player $1$ and $\tau$ be a strategy of
player $2$, then $\gamma_n(\sigma, \tau)$ is nondecreasing in $n$.
This implies that the corresponding sequence of values $(v_n)_{n\in
\mathbb{N}
}$ is nondecreasing in $n$. Denote $v=\sup_n v_n$ and let us show
that $v$ is the uniform value.

Fix $\varepsilon>0$, consider $N$ such that $v_N\geq v-\varepsilon$ and
$\sigma^*$ a strategy of player 1 which is optimal in $\Gamma_N$. We
have for each $\tau$ and, for every $n\geq N$,
\[
\gamma_n\bigl(\sigma^*, \tau\bigr) \geq\gamma_N\bigl(
\sigma^*, \tau\bigr)\geq v_N \geq v-\varepsilon.
\]
Hence, the strategy $\sigma^*$ guarantees $v-\varepsilon$ in the
uniform sense. This is true for every positive $\varepsilon$, thus player
$1$ guarantees $v$ in the uniform sense.

Using the monotone convergence theorem, we also have
\begin{eqnarray*}
\gamma^*\bigl(\sigma^*,\tau\bigr)&=&\mathbb{E}_{\sigma^*,\tau} \Biggl( \lim
_n \frac{1}{n}\sum_{t=1}^n
g(x_t,i_t,j_t) \Biggr) \\
&=& \lim
_n \mathbb{E}_{\sigma^*,\tau} \Biggl( \frac{1}{n}\sum
_{t=1}^n g(x_t,i_t,j_t)
\Biggr) \\
&\geq &  v-\varepsilon.
\end{eqnarray*}
We now show that player $2$ can also guarantee $v$ in the uniform
sense. Consider for every $n$, the set
\[
K_n=\bigl\{\tau, \forall\sigma,  \gamma_n(\sigma,
\tau)\leq v\bigr\}.
\]
$K_n$ is nonempty because it contains an optimal strategy for player 2
in $\Gamma_n$ (since $v_n\leq v$). The set of strategies of player 2 is
compact, hence by continuity of the $n$-stage payoff $\gamma_n$, $K_n$
is itself compact. $\gamma_n \leq\gamma_{n+1}$ implies
$K_{n+1}\subset
K_n$ hence $\bigcap_n K_n\neq\varnothing$: there exists $\tau^*$
such that for every strategy of player $1$, $\sigma$ and for every
positive integer $n$, $\gamma_n(\sigma,\tau) \leq v$.
It follows that both players can guarantee $v$, thus $v$ is the uniform
value.

By the monotone convergence theorem, we also have
%
\[
\gamma^*\bigl(\sigma,\tau^*\bigr)= \mathbb{E}_{\sigma,\tau^*} \Biggl( \lim
_n \frac{1}{n}\sum_{t=1}^n
g(x_t,i_t,j_t) \Biggr)
 = \lim
_n \mathbb{E}_{\sigma,\tau^*} \Biggl( \frac{1}{n}\sum
_{t=1}^n g(x_t,i_t,j_t)
\Biggr) \leq v.
\]
Hence, $v$ is the sup and limsup value.
\end{pf*}

\begin{remark}
The fact that the sequence of $n$-stage values $(v_n)_{n\geq1}$ is
nondecreasing is not enough to ensure the existence of the uniform
value. For example, consider the Big Match \cite{Blackwell68} with no
signals: $v_n=1/2$ for each $n$, but there is no uniform value.
\end{remark}

\begin{remark}
The theorem states the existence of a $0$-optimal strategy for
player 2 but player 1 may only have $\varepsilon$-optimal
strategies. For example, in the following MDP, there are two
absorbing states, two nonabsorbing states with payoff~$0$ and two
actions $\mathit{Top}$ and $\mathit{Bottom}$:
\[
\begin{tabular}{cc}
$\lleft(
\begin{array}{c}
1/2 (s_1) + 1/2 (s_2)\vspace*{1pt} \\
0^*\\
\end{array}
 \rright)$ & $\lleft(
\begin{array}{c}
s_2 \\
1^*\\
\end{array}
 \rright)$. \\
$s_1$ & $s_2$
\end{tabular}
\]
The starting state is $s_1$ and player $1$ observes nothing. A good
strategy is to play $Top$ for a long time and then $Bottom$. While
playing $Bottom$, the process absorbs and with a strictly positive
probability the absorption occurs in state $s_1$ with absorbing payoff
$0$. So
player $1$ has no strategy which guarantees the uniform value of 1.
\end{remark}






\printaddresses
\end{document}